\documentclass[12pt]{article}
\usepackage{latexsym}
\usepackage{enumerate}
\usepackage{epsf,epsfig,amsfonts,a4wide}
\usepackage{amsfonts}
\usepackage{url}
\usepackage{amsmath,amssymb,amsthm}
\usepackage{epsf,epsfig,amsfonts,graphicx}
\usepackage{a4wide}
\usepackage{lscape}

\usepackage[centerlast]{caption}
\usepackage{graphics,amsmath,amssymb,amscd}
\usepackage{amsbsy}
\usepackage{amsthm}
\usepackage{gensymb}
\usepackage{float}
\usepackage{graphicx}

\usepackage{hyperref}
\usepackage{epstopdf}
\usepackage{subfig}

\newcommand{\be}{\begin{equation}}
\newcommand{\ee}{\end{equation}}

\usepackage{amsthm}
\usepackage{mathdots}

\usepackage{amscd}

\newcommand{\mR}{\mathbb{R}}
\newcommand{\mS}{\mathbb{S}}

\newcommand{\om}{\omega}
\newcommand{\Om}{\Omega}
\newcommand{\eps}{\varepsilon}
\newcommand{\ph}{\varphi}

\newcommand{\me}{\mathrm{e}}
\newcommand{\mi}{\mathrm{i}}
\newcommand{\dif}{\mathrm{d}}

\title{On phenomenon of scattering on resonances associated with discretisation of systems with fast rotating phase}
\author{Anatoly Neishtadt$^{1,2}$\footnote{E-mail: {\tt A.Neishtadt@lboro.ac.uk}}, Tan Su$^1$\footnote{E-mail: {\tt T.Su@lboro.ac.uk}}\\
$^1$ Loughborough University, Loughborough, LE11 3TU, UK\\
 $^2$ Space Research Institute, Moscow, 113997, Russia}

\begin{document}
\maketitle

\begin{abstract}
Numerical integration of ODEs by standard numerical methods reduces a continuous time problems to discrete time problems. Discrete time problems have intrinsic properties that are absent in continuous time problems. As a result, numerical solution of an ODE may demonstrate dynamical phenomena that are absent in the original ODE. We show that numerical integration of system with one fast rotating phase lead to a situation of such kind: numerical solution demonstrate phenomenon of scattering on resonances that is absent in the original system. 
\vskip 10pt
\noindent MSC2010 numbers: {\tt 34F15, 34C29, 65D30}

\end{abstract}

\vskip 20pt

Numerical integration of ODEs reduces a continuous time problems to discrete time problems. For arbitrarily small time step of a numerical method such discrete time problem may have intrinsic properties that are absent in the original continuous time problem (see, e.g., \cite{Fiedler} and references therein). In this note we describe a situation of such kind that was not reported before: appearance of scattering on resonances in numerical integration of systems with one fast rotating phase. 


The system under consideration has the form 
\begin{eqnarray}\label{original1}
\begin{cases}
\dot I=f(I, \ph, \eps)\,, \quad I\in U\subseteq \mR^l \,,\\
\dot\ph=\dfrac{\omega(I)}{\eps}+g(I, \ph, \eps)\,, \quad \ph\in \mS^1\ {\text {mod}}\ 2\pi.
\end{cases}
\end{eqnarray}
Here $\ph$ is an angular variable (phase), functions $f$ and $g$ are $2\pi$-periodic in $\ph$, $\eps$ is a small positive parameter. We assume that values of function $\om$ are separated from $0$ in the domain $U$ by a positive constant, and that functions $\om, f,g$ are real-analytic in $U\times\mS^1\times [0,\eps_0], \, \eps_0={\text {const}}>0$. In system (\ref{original1}) the variables $I,\ph$ are called slow and fast variables, respectively, the function $\om/\eps$ is called a frequency. System (\ref{original1}) is called a system with fast rotating phase \cite{Bogol}. Dynamics in many problems is described by systems of form (\ref{original1}).

Evolution of slow variables $I$ in system (\ref{original1}) on time intervals of lengths of order 1 is approximately described by the averaged system:
\begin{equation}\label{aver}
\dot J = F(J),\quad F(J)=\frac{1}{2\pi}\int\limits_0^{2\pi}f(J,\ph,0)\,\dif\ph 
\end{equation}
(the averaging method, see, e.g., \cite{Bogol}). One can use higher approximations of the averaging method as well \cite{Bogol}. We now discuss a phenomenon that appears if system (\ref{original1}) is solved numerically. 

Numerical integration of system (\ref{original1}) with a fixed time step $\kappa$ effectively introduces into dynamics a new ``numerical'' frequency $\Om=2\pi/\kappa$. 
In the process of evolution of slow variables $I$ the frequency $\omega(I)/\eps$ changes and passes through resonances with the numerical frequency $2\pi/\kappa$. Passage through a resonance leads to a scattering on resonance: a deviation of dynamics of $I$ from the averaged dynamics described by (\ref{aver}) (see, e.g., \cite{scat}, a scattering on resonance for the first time was discussed in \cite{Chir}). This deviation depends on value of phase $\ph$ at the moment of passage through resonance. Consider for example, the simplest case when the system is integrated by the Euler method. Then the values of $I, \ph$ at the moments of time $n\kappa$ and $(n+1)\kappa$ are related as
\begin{equation} \label{E_map}
I_{n+1}=I_n+\kappa f(I_n, \ph_n, \eps),\quad \ph_{n+1}=\ph_n+\kappa\frac{\omega(I_n)}{\eps}+\kappa g(I_n, \ph_n, \eps).
\end{equation}
We can say that this numerical procedure integrate the time-periodic system of ODEs
$$\dot I=\kappa f(I, \ph, \eps)\sum\limits_{n=-\infty}^{\infty}\delta(t-n\kappa),\quad \dot\ph=\kappa \left(\frac{\omega(I)}{\eps}+g(I, \ph, \eps)\right)\sum\limits_{n=-\infty}^{\infty}\delta(t-n\kappa),
$$
where $\delta(\cdot)$ is $\delta$-function. We can use the standard identity 
$$ \sum\limits_{n=-\infty}^{\infty}\delta(t-n\kappa) =\frac{1}{\kappa}\sum\limits_{n=-\infty}^{\infty}\me^{\mi2\pi nt/\kappa}.
$$
Consider Fourier series expansion for $f$:
$$f(I,\ph,0)=\sum\limits_{n=-\infty}^{\infty}f_n(I)\me^{\mi\ph}\,\,.
$$
Near a low order resonance $n_1\frac{\omega(I)}{\eps}+n_2\frac{2\pi}{\kappa}=0$, where $n_1$ and $n_2$ are co-prime integer numbers, the dynamics is approximately described by the partially averaged equations (see \cite{scat})
\begin{equation}\label{paver}
\dot I =F(I)+\sum\limits_{\substack{q=-\infty\\ q\ne 0}}^{\infty} f_{qn_1}(I)\,\me^{\mi q(n_1\ph+n_2\Om t)},\quad \dot\ph =\frac{\omega(I)}{\eps}\, .
\end{equation}
For such form of equations there are asymptotic formulas for amplitude of scattering on resonance in different situations (see \cite{Kev, scat}). 
Assume that $f$ is analytic in a strip $|{\text {Im}} \,\ph|\le \sigma$. Then $|f_n| < {\text {const}}\cdot\me^{-\sigma |n|} $. Then under rather general assumptions the amplitude of the scattering is $\sim \sqrt\eps\, \me^{-\sigma |n_1| }$. For $n_2=\pm1$ we get that the amplitude of the scattering is $\sim \sqrt\eps\, \me^{-\sigma \frac{2\pi\eps}{\omega(I)\kappa} }$. Thus the considered effect is exponentially small in natural circumstances when $\kappa\ll\eps$.

We first will demonstrate the existence of scattering on resonance for iterations of maps of the form (\ref{E_map}). To this end we will integrate numerically by the Euler method with unrealistically big fixed time step the following simple system of equations with one rotating phase:
\begin{equation}\label{test1}
\dot I_1=1,\quad \dot I_2 =\cos\ph,\quad \dot \ph =\frac{\omega(I_1)}\eps.
\end{equation} 
If this system is integrated by the Euler method, then, according to formula (\ref{paver}), near the resonance $\frac{\om(I_1)}\eps+n_2\Om=0$, the dynamics is approximately described by the system
\begin{equation}
\dot I_1=1,\quad \dot I_2=\cos(\ph+n_2\Om t),\quad \dot\ph=\frac{\om(I_1)}\eps.
\end{equation}
Let the resonance take place at $t=t_*$, i.e. $\frac{\om(I_1(t_*))}\eps+n_2\Om=0$. Denote $I_{1*}=I_1(t_*)$, $\ph_*=\ph(t_*)$. Then passage through a narrow neighbourhood of this resonance produces a jump in $I_2$ given by the asymptotic formula
\begin{equation}\label{jump}
\Delta I_2 =\sqrt{\frac{2\pi\eps}{|\om_*'|}}\cos\!\left(\ph_*+n_2\Om t_*+{\rm sgn}(\om_*') \tfrac\pi4\right)+O(\eps^{\frac32})
\end{equation}
where $\om_*'= \om'(I_{1*})$ (c.f. \cite{Bosley} for the estimate of the error term). 
 Here $\om'={\partial \om}/{\partial I_1}$, we assume that $\om'(I_{1*})\ne0$. This phenomenon is called a ``scattering'' because the value of jump depends on the value of phase $\ph=\ph_*$ at the moment of passage through a resonance. For discrete time systems these jumps were first observed in \cite{PF}.

\begin{figure}
\centering
\subfloat[$\kappa=\eps$]{\label{euler1}\includegraphics[height=1.8in,width=5in]{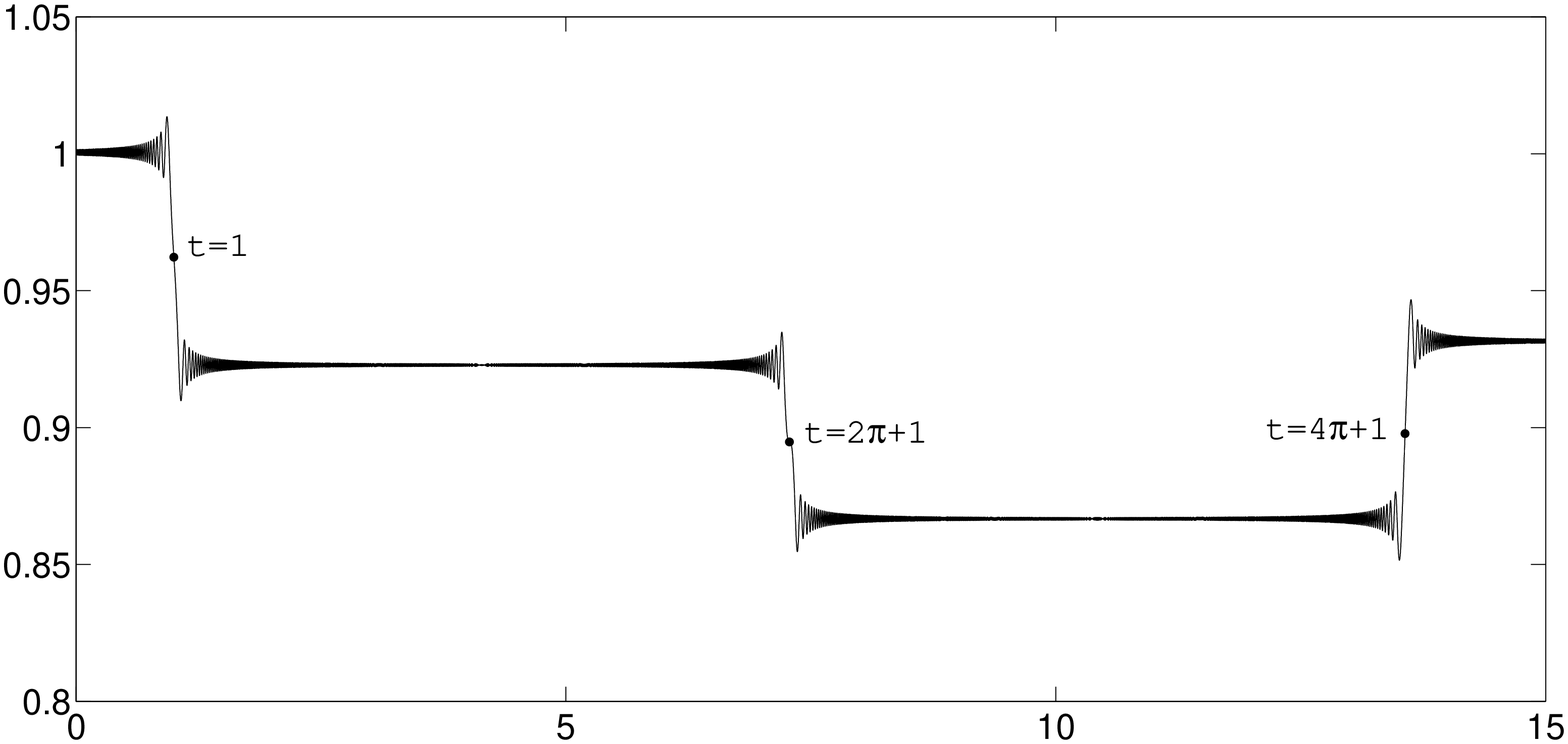}} \\ 
\subfloat[$\kappa=\eps/2$]{\label{euler2}\includegraphics[height= 1.8in,width=5in]{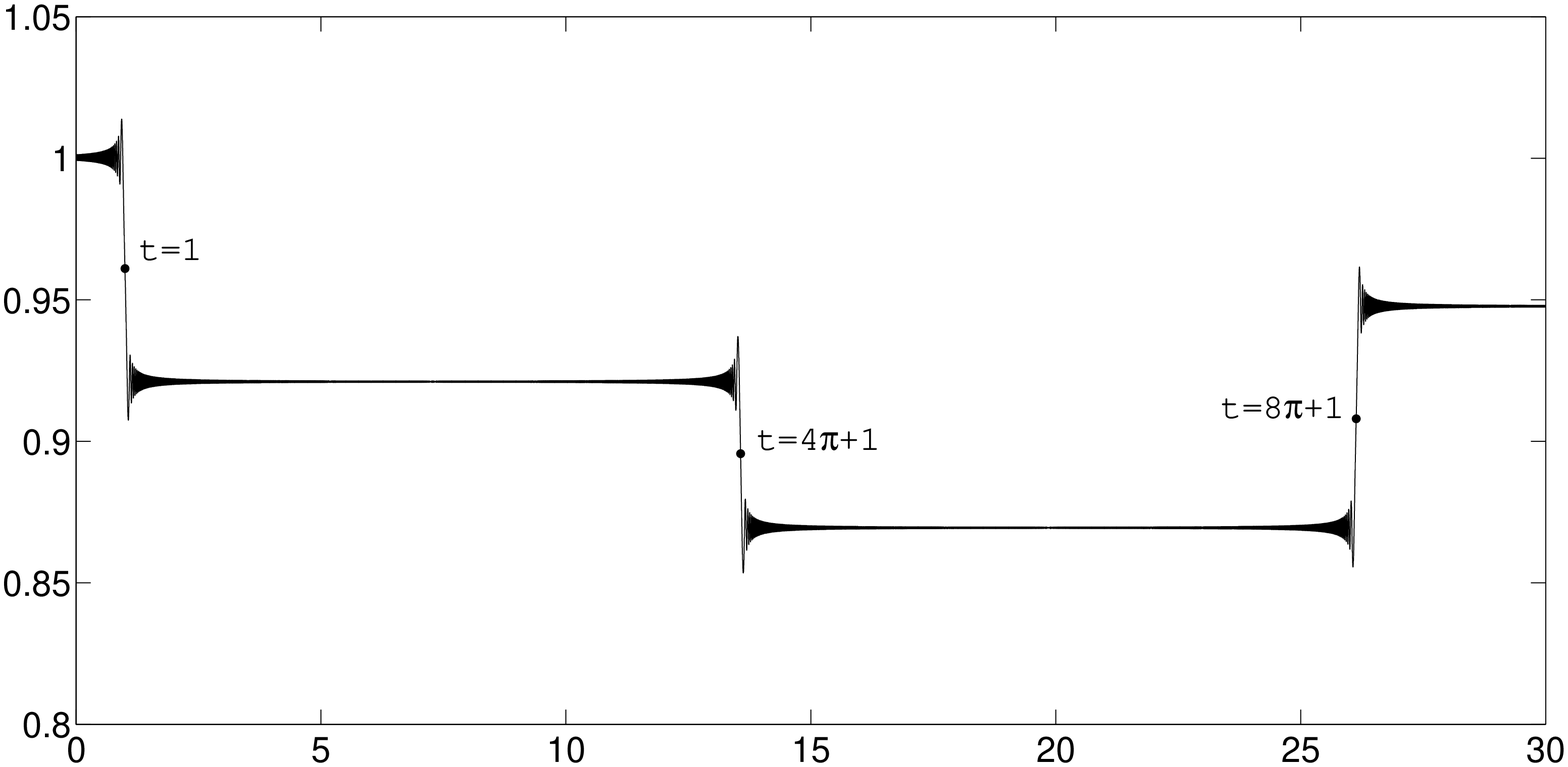}} \\ 
\subfloat[$\kappa=\eps/5$]{\label{euler5}\includegraphics[height= 1.8in,width=5in]{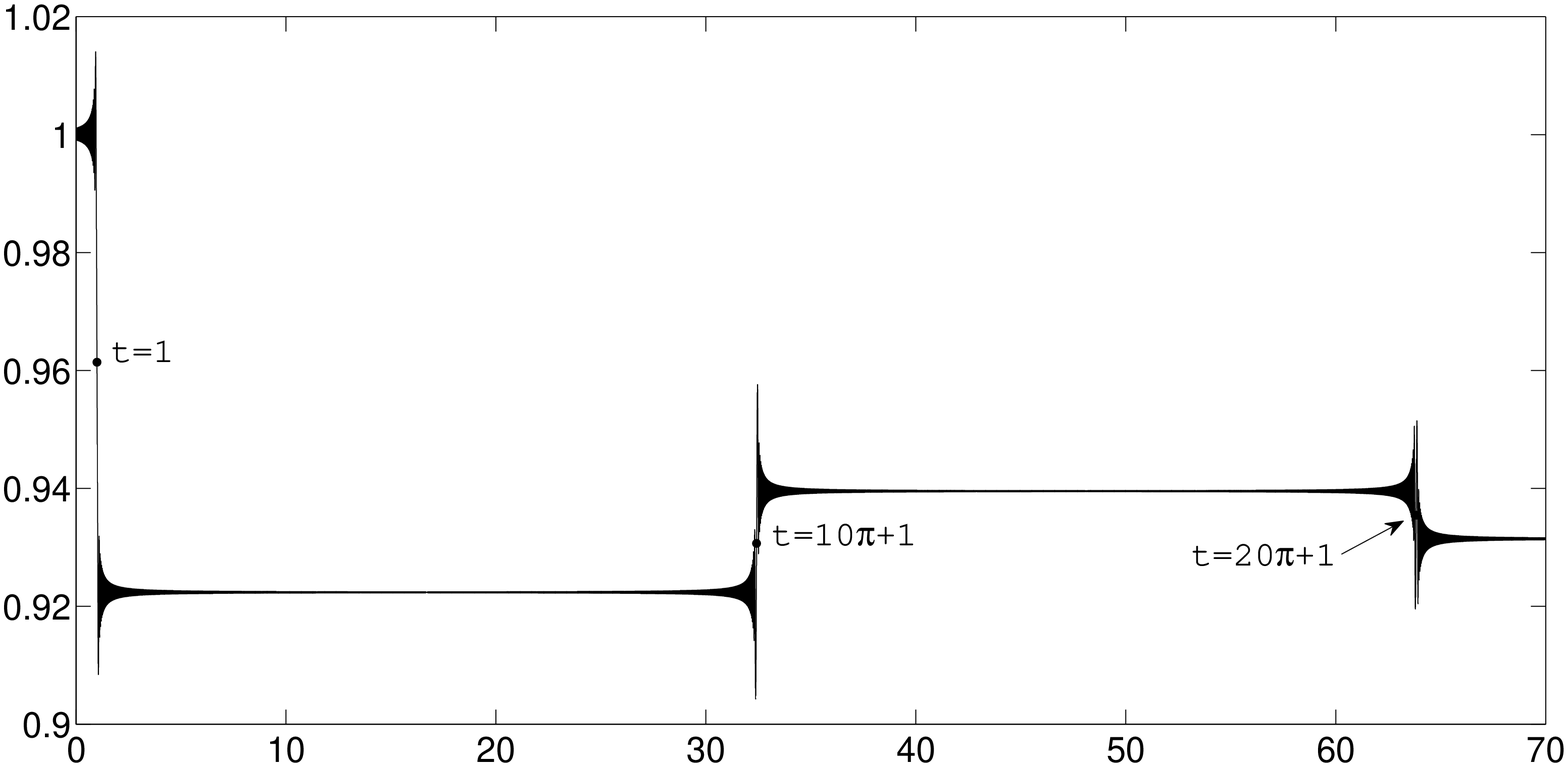}} \\ 
\subfloat[$\kappa=\eps/10$]{\label{euler10}\includegraphics[height= 1.8in,width=5in]{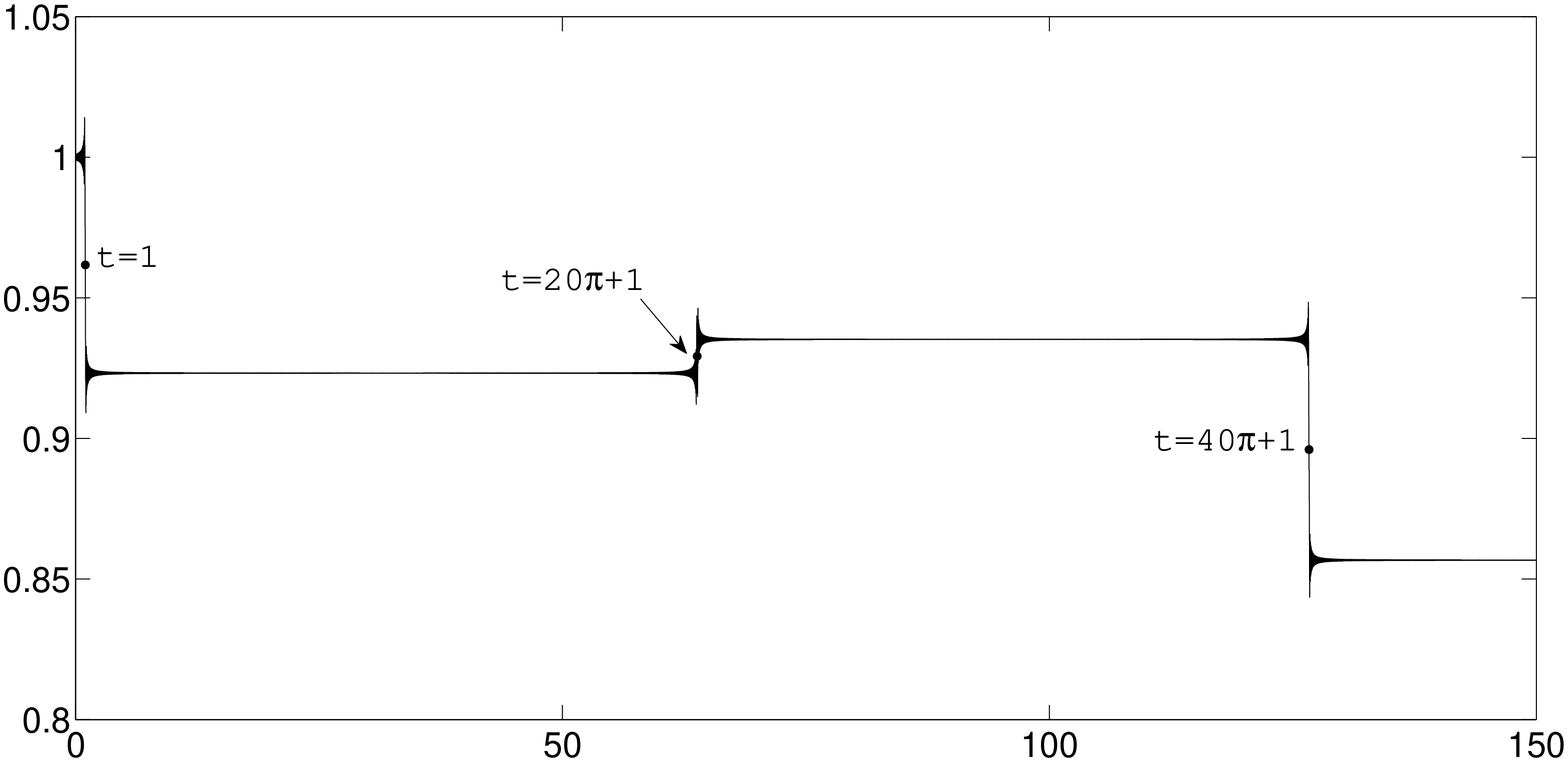}} 
\caption{Euler discretisation, $\dot I_2=\cos\ph$, $\om=t-1$, $\eps=0.001$}
\label{euler}
\end{figure}

We take $\omega=I_1$ and the following values of the time step $\kappa$: $\eps$, $\eps/2$, $\eps/5$, $\eps/10$. Take $\eps=0.001$.
 Initial conditions for each run are $I_1(0)=-1$, $I_2(0)=1$, $\ph(0)=0$. Thus $\omega=I_1=t-1$, and at $t=1$ the system passes through the actual resonance $\omega=0$, where $I_2$ has a jump. Numerical results are presented at Fig. \ref{euler}. One can see a jump in $I_2$ at $t=1$ due to the actual resonance, as well as jumps at $t=1+2\pi n (\eps/\kappa), \,\, n=1,2,$ due to the discretisation. (Plots here and in all other figures are done with MATLAB.) 
 




Comparison of calculation of jumps using formula (\ref{jump}) with numerical results in Fig. \ref{euler} is presented in Table 1 \footnote{ 
Value $\ph_*$ in Table \ref{table1} was calculated as follows: $\ph_*=\ph_1+\eps^{-1}\int_{t_1}^{t_*}\om\,\dif t =\ph_1+\eps^{-1}\left((t_*-1)^2-(t_1-1)^2\right)/{2}$, where $t_1$ is the time moment in the considered discrete grid immediately preceding to $t_*$, and $\ph_1$ is the value of $\ph$ for the numerical solution at $t=t_1$.}.

\begin{table}[H]
\centering
\begin{tabular}{|l|l|l|l|l|}
\hline
$\ \kappa$ & $\ t_*$ & $\ \ph_*$ & $\ \Delta I_2$ (theoretical) & $\ \Delta I_2$ (numerical)\\
\hline
$\ \eps$ & $\ 4\pi+1$ & $\ 78462.611$ & $\ 0.0653$ & $\ 0.0649$\\
\hline
$\ \eps/2$ & $\ 8\pi+1$ & $\ 315333.34308$ & $\ 0.0786$ & $\ 0.0782$\\
\hline
$\ \eps/5$ & $\ 20\pi+1$ & $\ 1973427.065167$ & $\ -0.0082$ & $\ -0.0080$\\
\hline
$\ \eps/10\ $ & $\ 40\pi+1\ $ & $\ 7895189.742154\ $ & $\ -0.0784$ & $\ -0.0785$\\
\hline
\end{tabular}
\caption{Theoretical and numerical results for jump, Euler discretisation}
\label{table1}
\end{table}

\begin{figure}
\centering
\subfloat[$\kappa=\eps$]{\label{RKcos1}\includegraphics[height= 1.8in,width=5in]{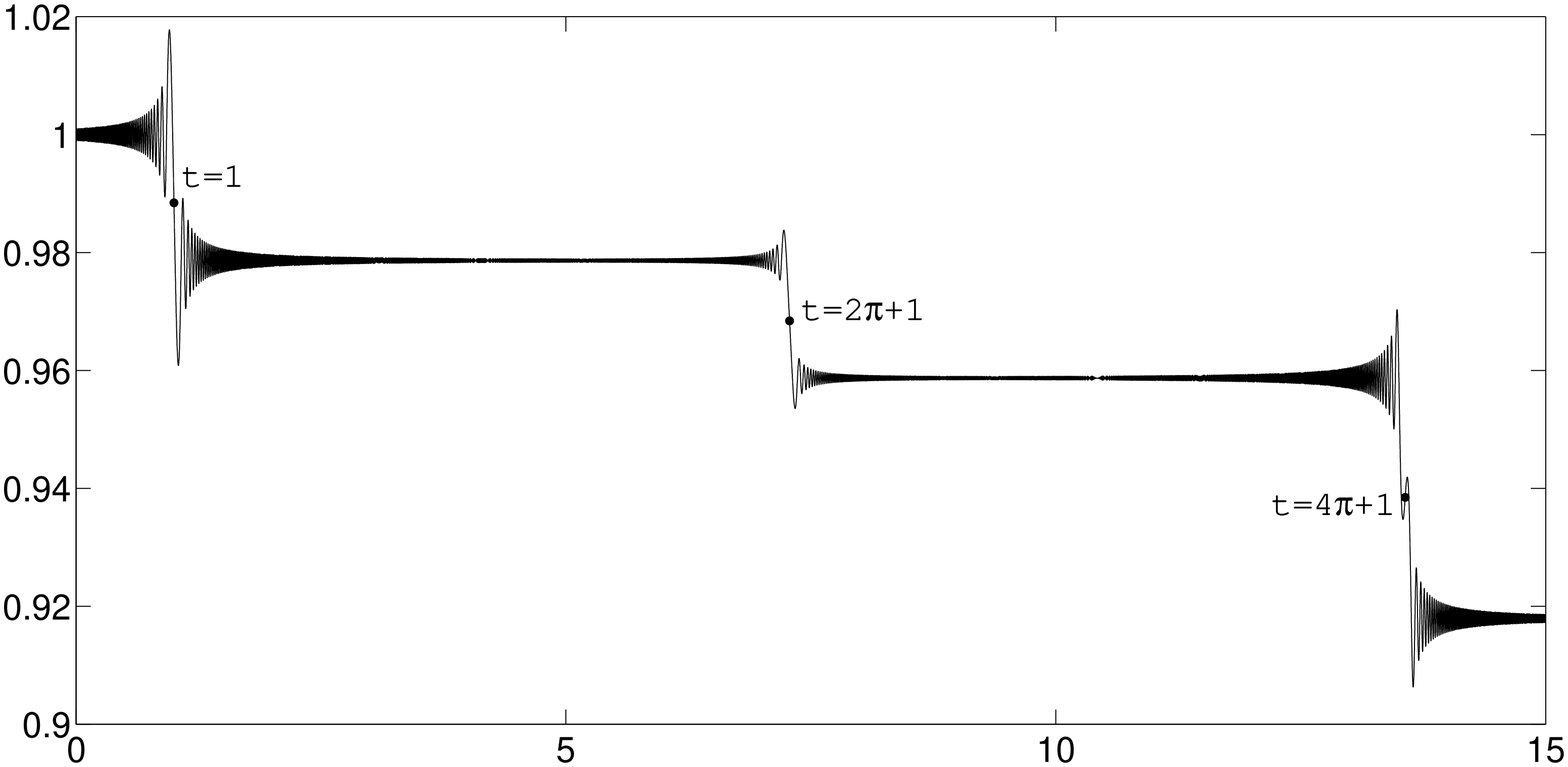}} \\ 
\subfloat[$\kappa=\eps/2$]{\label{RKcos2}\includegraphics[height= 1.8in,width=5in]{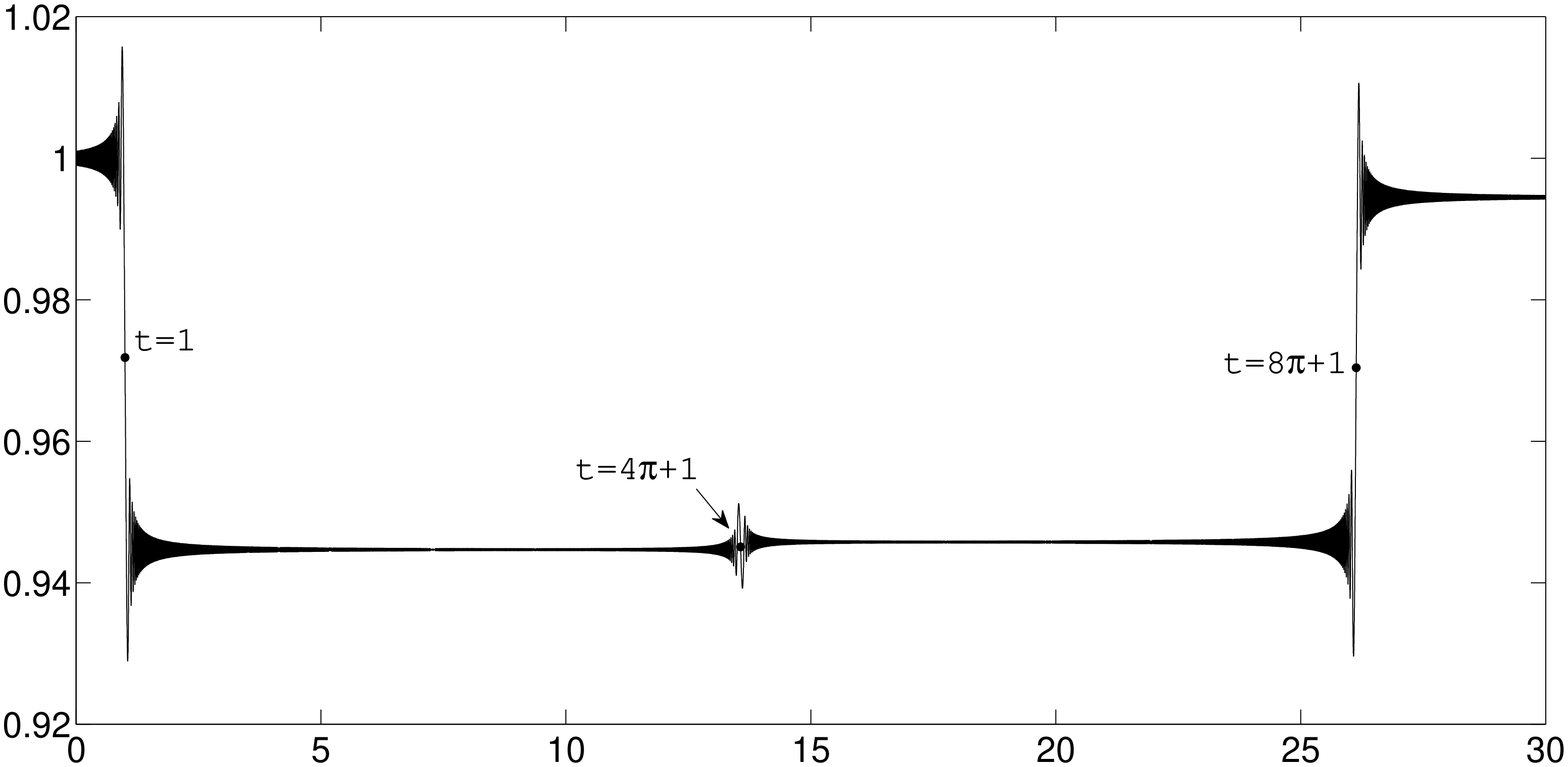}} \\ 
\subfloat[$\kappa=\eps/5$]{\label{RKcos5}\includegraphics[height= 1.8in,width=5in]{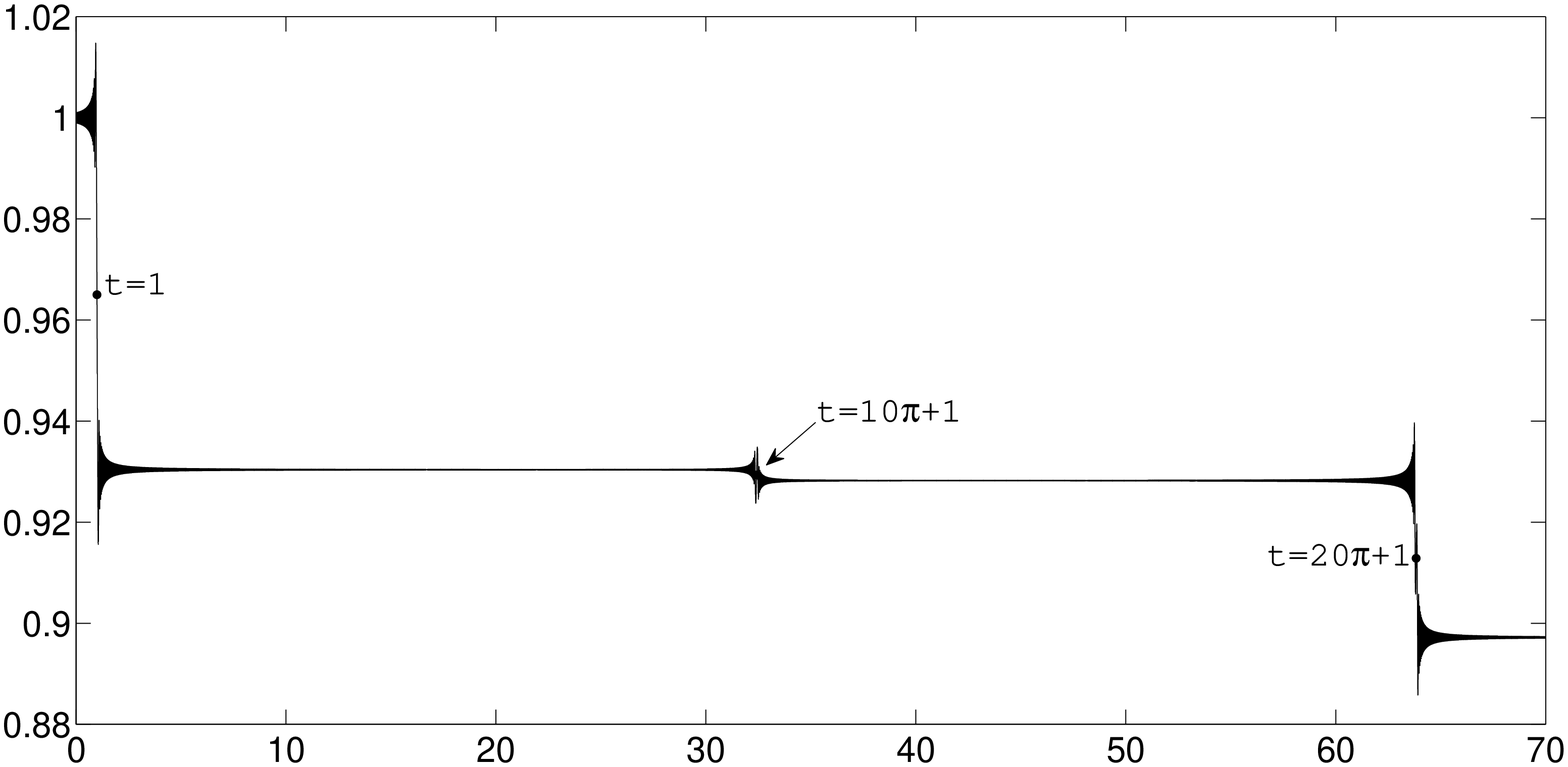}} \\ 
\subfloat[$\kappa=\eps/10$]{\label{RKcos10}\includegraphics[height= 1.8in,width=5in]{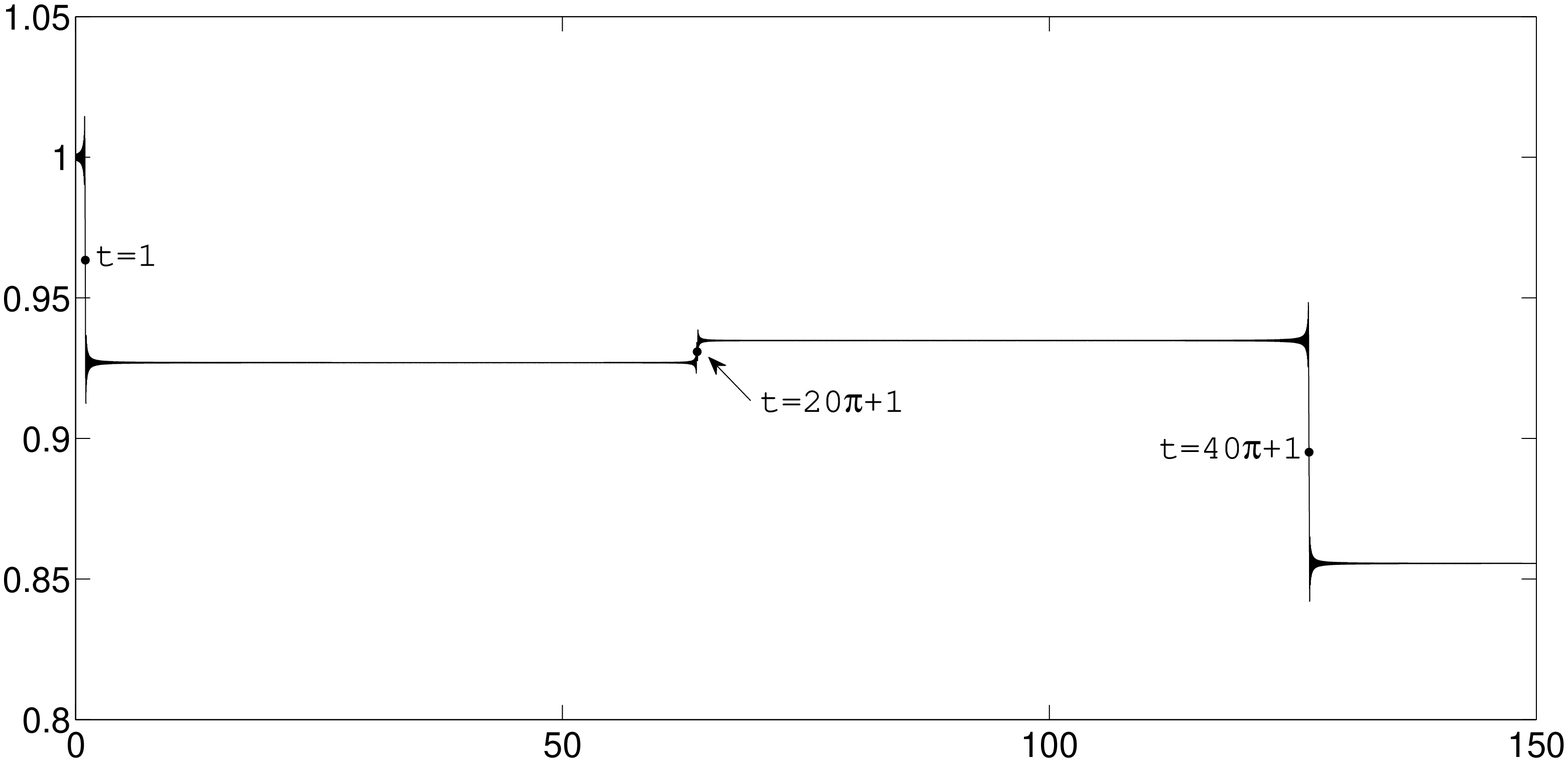}} 
\caption{Runge-Kutta discretisation, $\dot I_2=\cos\ph$, $\om=t-1$, $\eps=0.001$}
\label{RKcos}
\end{figure}

Similar scatterings on resonances occur for other maps with fast rotating phase. We demonstrate this for the map that appears in solving system (\ref{test1}) by 4th order Runge-Kutta method. We again take $I_1(0)=-1,\, I_2(0)=1, \, \ph(0)=0,\, \eps= 0.001$ and the same unrealistically big values of the time step $\kappa$: $\eps$, $\eps/2$, $\eps/5$, $\eps/10$. Results are presented in Fig. \ref{RKcos}. One can see jumps in $I_2$ at the the same moments of time as in Fig. \ref{euler}.





\begin{figure}
\centering
\includegraphics[width=16cm]{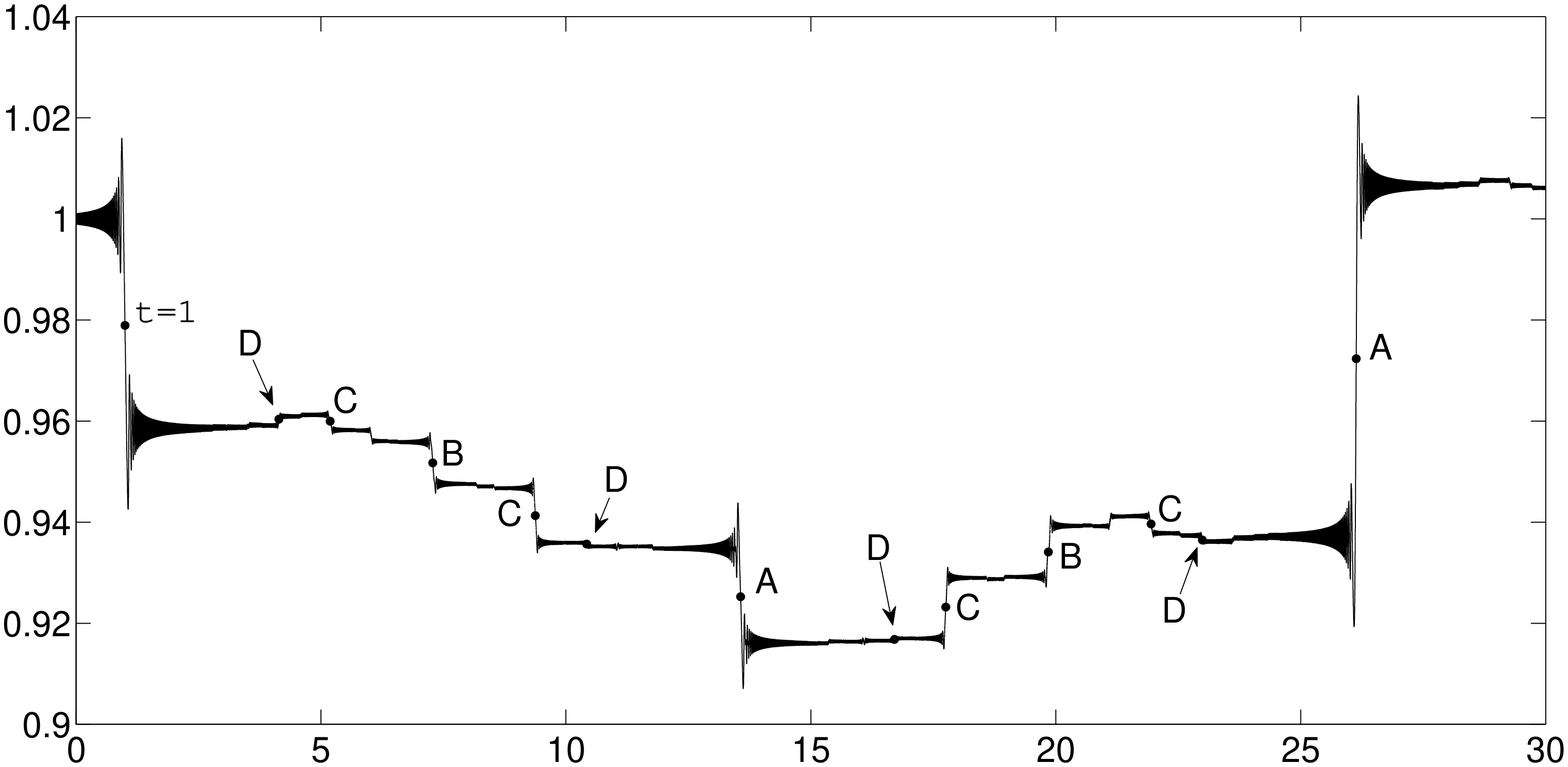}
\caption{$\dot I_2 =\cos\ph+\frac12\cos2\ph+\cdots+\frac1{2^{10}}\cos11\ph$, $\om=t-1$, $\eps=0.001$, $\kappa=\eps/2$}
\label{RKcos11}
\end{figure}

Now let us consider scatterings on resonances for solving the following system by 4th order Runge-Kutta method:
\begin{equation}\label{test11}
\dot I_1=1, \quad \dot I_2 =\cos\ph+\frac12\cos2\ph+\cdots+\frac1{2^{10}}\cos11\ph, \quad \dot\ph=\frac{\om(I_1)}{\eps}
\end{equation}
We again take $ I_1(0)=-1,\, I_2(0)=1, \, \ph(0)=0,\, \eps= 0.001$. Result for $\kappa=\eps/2$ is shown in Fig. \ref{RKcos11}. We use this unrealistically big value of the time step to demonstrate existence of the phenomenon. 

Resonances should appear at $n_1\frac{\om}\eps+n_2\frac{2\pi}\kappa=0$, where $n_1$ and $n_2$ are co-prime integer numbers. Using $\om=t-1$, $\kappa=\eps/2$, we obtain $t=-\frac{n_2}{n_1}\cdot4\pi+1$ at resonances.

We analyse the right hand side of $\dot I_2$ term by term. For the first term, $\cos\ph$, we have $n_1=1$, the resonances occur at the points $t=k\cdot4\pi+1$, $k=1,2,3,\ldots$ (marked as A in Fig. \ref{RKcos11}). For the next term, $\frac12\cos2\ph$, we have $n_1=2$. The possible values of $n_2$ are $n_2=-1,-3,-5,\ldots$ and the time moments for resonances are $t=(2k-1)\cdot2\pi+1$, $k=1,2,\ldots$ (shown in Fig. \ref{RKcos11} as B, which are $2\pi+1$, $4\pi+1$, \ldots). Similarly, for the term $\frac14\cos3\ph$, $n_1=3$, $n_2=-1,-2,-4,-5,\ldots$, $t=\frac43\pi+1,\frac83\pi+1,\frac{16}3\pi+1, \frac{20}3\pi+1, \ldots$. For the term $\frac18\cos4\pi$, $n_1=4$, $n_2=-1,-3,-5,-7,\ldots$, and resonances are at $t=\pi+1, 3\pi+1, 5\pi+1, 7\pi+1, \ldots$. The corresponding points are marked in Fig. \ref{RKcos11} as C and D, respectively. If we go on with this procedure, we will find every resonance for 11 terms by increasing $n_1$ from 1 to 11.

\begin{figure}
\centering
\subfloat[$\kappa=\eps/2$]{\label{RKcosinf2}\includegraphics[height= 1.8in,width=5in]{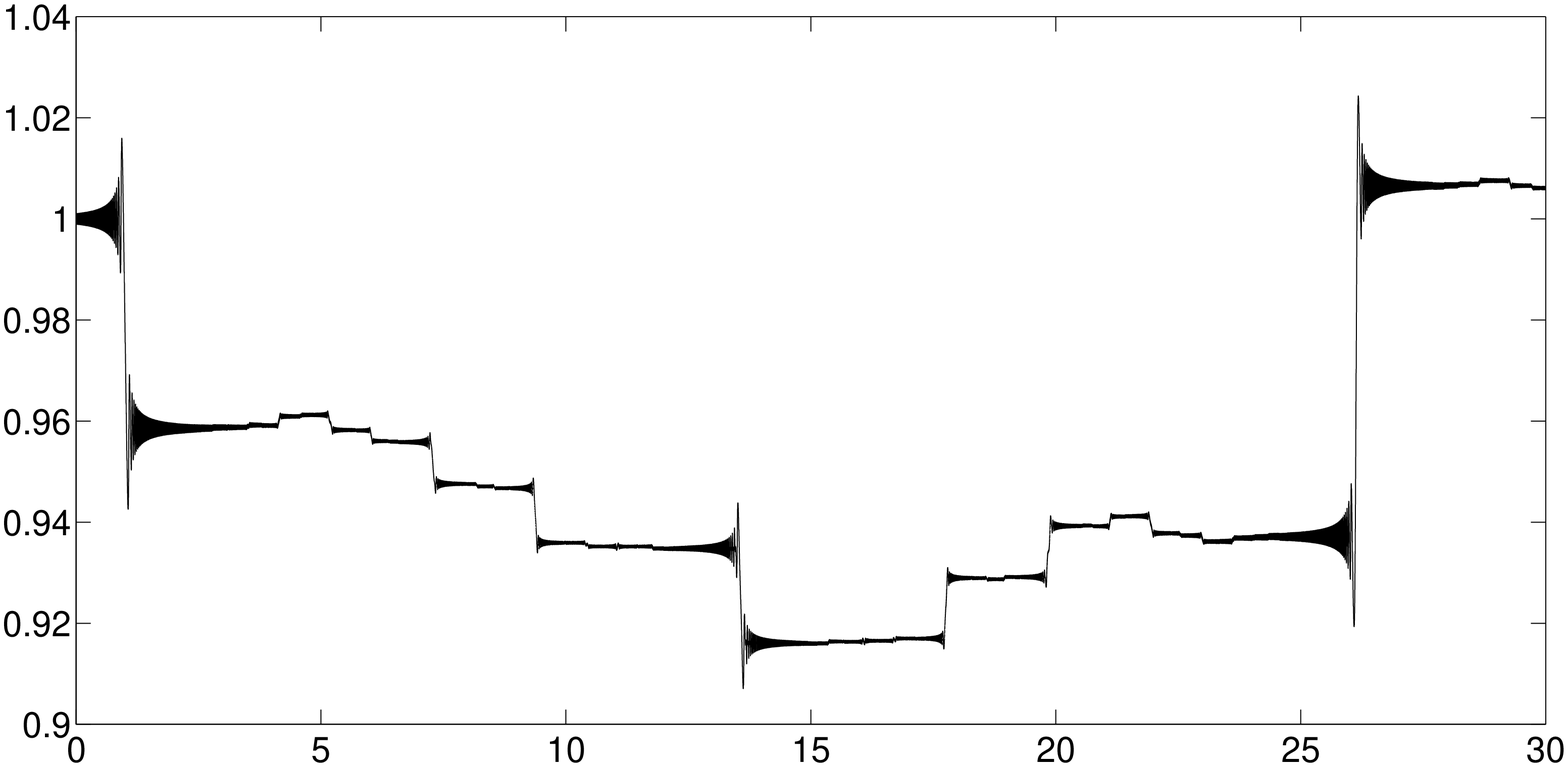}} \\ 
\subfloat[$\kappa=\eps/10$]{\label{RKcosinf10}\includegraphics[height= 1.8in,width=5in]{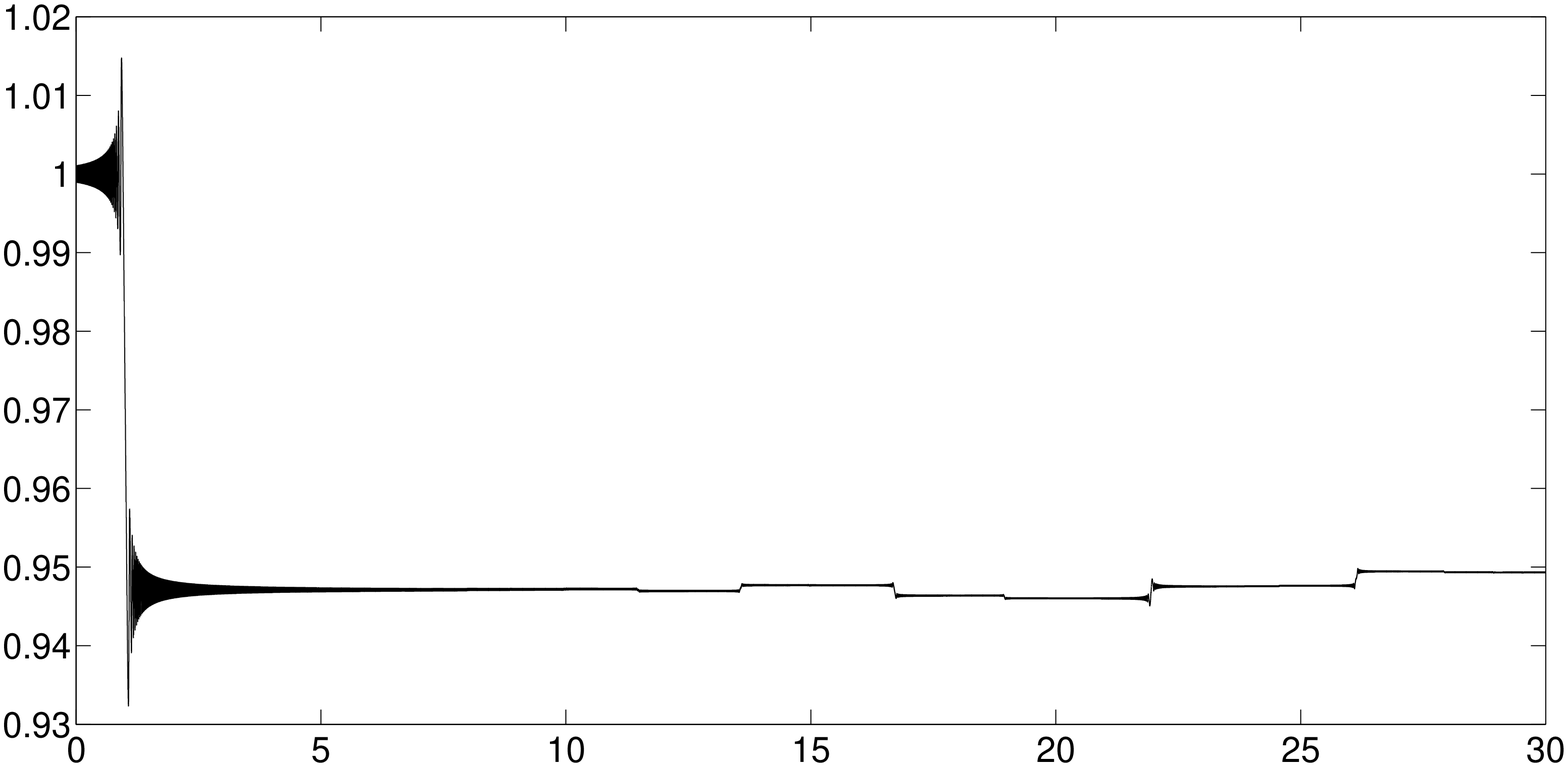}} \\ 
\subfloat[$\kappa=\eps/20$]{\label{RKcosinf20}\includegraphics[height= 1.8in,width=5in]{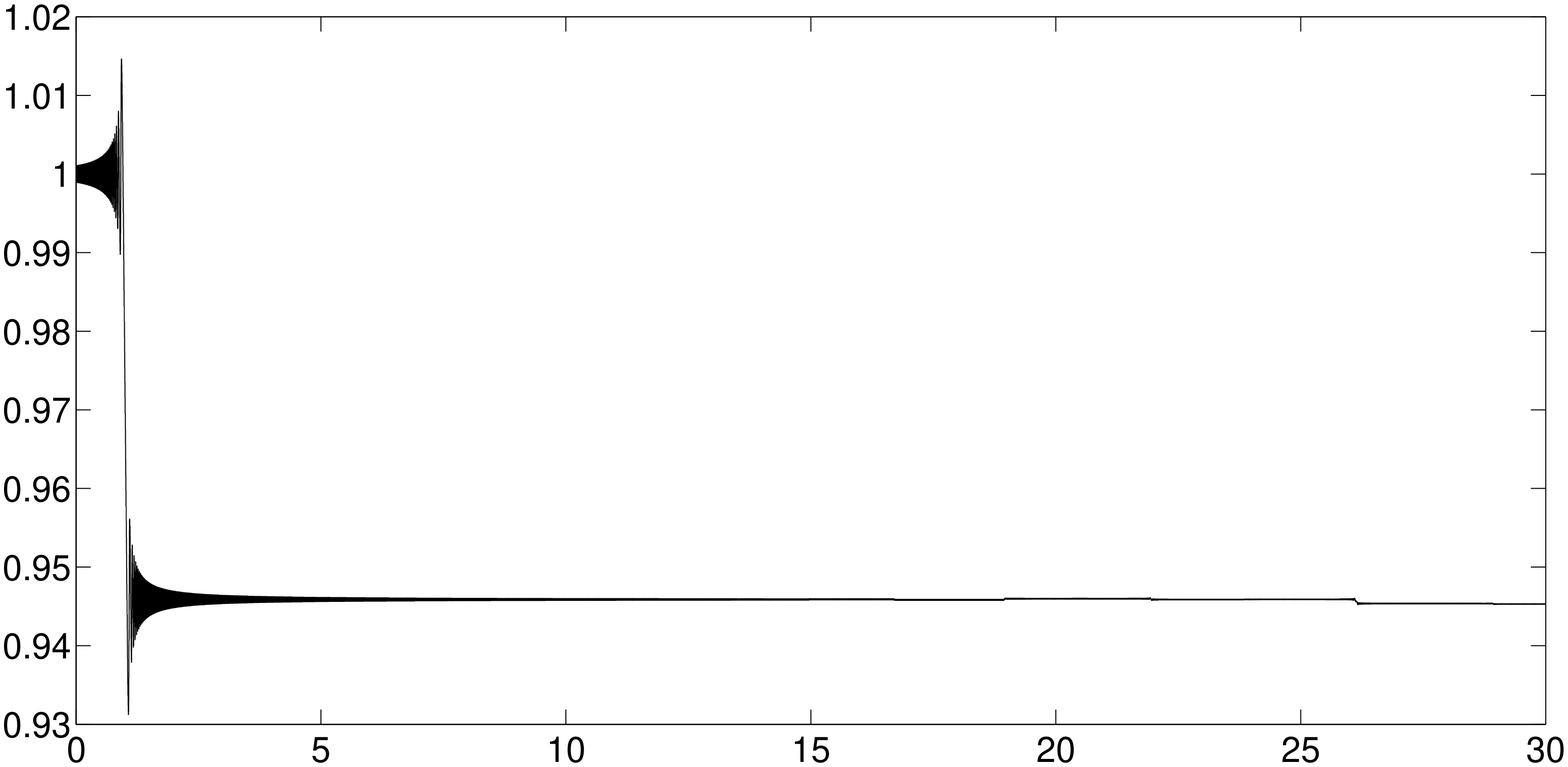}} \\
\caption{$\dot I_2=\dfrac{4\cos\ph-2}{5-4\cos\ph}$, $\om=t-1$, $\eps=0.001$, $\kappa=\eps/2$}
\label{RKcosinf}
\end{figure}

Now let us consider the case 
\begin{equation}\label{testinf}
\dot I_1=1, \quad \dot I_2 =\cos\ph+\frac12\cos2\ph+\frac14\cos3\ph+\cdots, \quad \dot\ph=\frac{\om(I_1)}{\eps}
\end{equation}
with infinite number of terms in $\dot I_2$. Let us represent $\dot I_2$ as
\begin{eqnarray*}
\dot I_2 =\ \rm{Re}\,\left(\me^{\mi\ph}+\frac12\me^{2\mi\ph}+\frac14\me^{3\mi\ph}+\cdots\right)
=\rm{Re}\,\left(\frac{\me^{\mi\ph}}{1-\frac12\me^{\mi\ph}}\right) =\frac{4\cos\ph-2}{5-4\cos\ph}.
\end{eqnarray*}
The result of numerical simulation by 4th order Runge-Kutta method for the same initial data and parameters as for system (\ref{test11}) is shown in Fig. \ref{RKcosinf2}. It is clearly seen that the resonances where the jump occurs are the same as for system (\ref{test11}) for the first four terms in $\dot I_2$. For later terms $\frac1{2^{n_1-1}}\cos n_1\ph$ in $\dot I_2$ we can find corresponding resonances, but the jumps become smaller and smaller due to decreasing of coefficients $\frac1{2^{n_1-1}}$ as $n_1$ increases. 

\begin{figure}
\centering
\includegraphics[width=16cm]{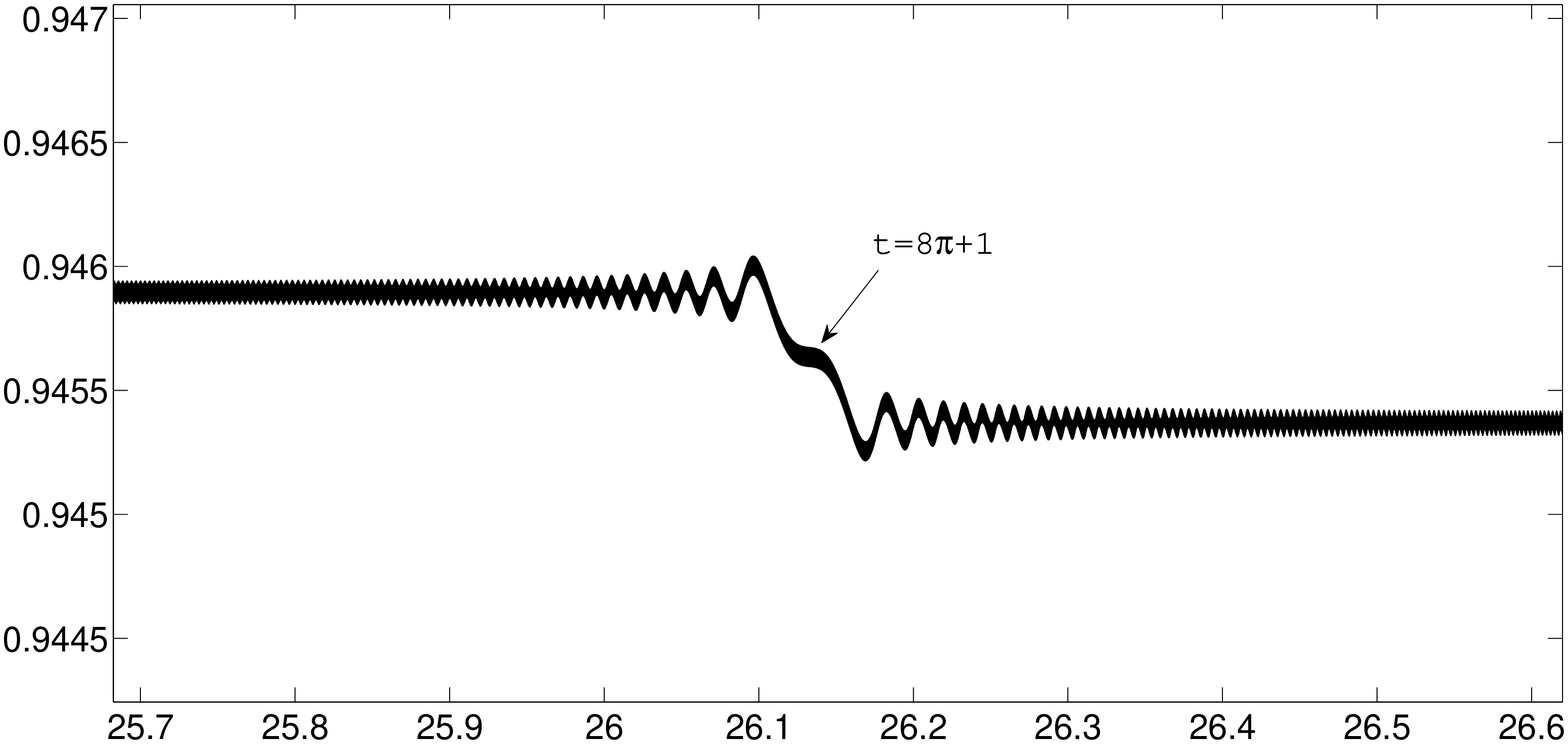}
\caption{$\dot I_2=\dfrac{4\cos\ph-2}{5-4\cos\ph}$, $\om=t-1$, $\eps=0.001$, $\kappa=\eps/20$, zoom in at $t=8\pi+1$}
\label{RKcosinfzoom}
\end{figure}

Finally, let us numerically integrate system (\ref{testinf}) by 4th order Runge-Kutta method with the same initial conditions and parameters, as in the previous run, but with the time step $\kappa=\eps/10$ and $\kappa=\eps/20$. The results are shown in Figs. \ref{RKcosinf10} and \ref{RKcosinf20}. One can see how amplitudes of jumps decay as $\kappa$ decays. We present a zoom of the jump at $t=8\pi+1$ in Fig. \ref{RKcosinfzoom}. This jump corresponds to the resonance $5\frac{\om}\eps-\frac{2\pi}\kappa=0$.
The curve in Fig. \ref{RKcosinfzoom} looks ``fat''. The reason is that the solution is the sum of a jump curve corresponding to the term $\cos5\ph$ in the right hand side of equation (\ref{testinf}), and high frequency oscillations corresponding to, mainly, terms $\cos m\ph,\,\,m=1,2,3,4\,,$ in (\ref{testinf}). If we magnify Fig. \ref{RKcosinfzoom} we would see these oscillations inside ``fat'' curve.

\medskip 
When behaviour of a system with fast rotating phase should be studied numerically, it looks as the most appropriate way is to use the averaging method and its higher approximations, like in \cite{Sanz-Serna}. If, however, the direct numerical integration of the original system with a standard constant step numerical integrator is used (e.g., to compare results with results obtained by the averaging), then the deviations of the numerical solution from the exact one have form of jumps on resonances between the internal frequency of the system and the frequency of discretisation.

\bigskip 

\end{document}